\documentclass[letterpaper, 10 pt]{article}  
\usepackage{times}
\usepackage{graphicx, amsmath, theorem, amssymb, algorithmic}

\newtheorem{theo}{Theorem}[section]

\newtheorem{rem}{Remark}[section]
\newtheorem{deff}{Definition}[section]

\def\be{\begin{equation}}
  \def\eeq{\end{equation}}
\def\<{\langle}
\def\>{\rangle}

\def\bold(#1){\textbf{#1}}
\def\roma(#1){\textrm{#1}}
\def\und(#1){\underline{#1}}
\def\ove(#1){\overline{#1}}
\def\mbold(#1){\mathbb{#1}}

\title{\LARGE \bf
Fault-Tolerant Control Allocation: an Unknown Input Observer based approach}

\author{Andrea Cristofaro and Tor Arne Johansen
\thanks{A. Cristofaro is with School of Science and Technology, University of Camerino, Italy. T.A. Johansen is with Department of Engineering Cybernetics, Center for Autonomous Marine Operations and Systems, Norwegian University of Science and Technology, Norway. email: {\tt andrea.cristofaro@unicam.it,  tor.arne.johansen@itk.ntnu.no}}%
}
\date{}

\begin{document}

\maketitle
\thispagestyle{empty}
\pagestyle{empty}


\begin{abstract}                          
This paper focuses on the use of unknown input observers for detection and isolation of actuator and effector faults with control reconfiguration in overactuated systems. The control allocation actively uses input redundancy in order to make relevant faults observable. The case study of an overactuated  marine vessel supports theoretical developments.
\end{abstract}

\section{Introduction}
The main objective of control allocation is to determine how to generate a specified control effect from a redundant set of actuators and effectors. Control effectors are devices or surfaces producing forces and moments, such as thrusters, propellers, fins or rudders, while actuators are electromechanical devices responsible to tune the magnitude, position and orientation of  single effectors. Due to input redundancy, several configurations leading to the same generalized force are admissible and for this reason the control allocation scheme commonly incorporates additional secondary objectives \cite{Bodson:02} \cite{FoJo:06} \cite{JoFo:13}, such as power or fuel consumption minimization. On the other hand, usually there are also some limitation factors to be accounted for: actuators/effectors dynamics, input saturation and other physical or operational constraints. One further advantage of actuator and effector redundancy is the possibility to reconfigure the control in order to cope with unexpected changes on the system dynamics, such as failures or malfunctions: in particular if the set of actuators and effectors is partially affected by faults, one can modify the control allocation scheme by preventing the use of inefficient/ineffective devices in the generation of control effect or compensating for the loss of efficiency. However, one key point for successfully re-allocating the control is the availability of adequate  information about the faults that have occurred; indeed, some accurate fault estimation and/or a correct identification of the faulty actuators or effectors are necessary to address the reconfiguration. Recent results toward fault tolerant control allocation are based on sliding-mode techniques   \cite{AlEd:06} \cite{Co-etal:05} and adaptive control strategies \cite{CaGa:10} \cite{TjJo:08}. 
Further investigations on this topic, with a more application-oriented character, are proposed for reconfiguration in flight control \cite{Bu-etal:01}, \cite{Zhang-etal:07} and fault accommodation in automated underwater vehicles \cite{Sarkar-etal:02}. An interesting bibliographical survey on the general problem of fault-tolerant control reconfiguration is provided in \cite{ZhJi:08}. \\
The aim of this paper is to present the use of unknown input observers for fault detection/isolation and control reconfiguration in overactuated systems. Unknown input observers \cite{CPZ:96} are a very useful tool for generating robust detection filters, as they can be made insensitive to certain input space directions if some structural algebraic conditions on the system are fulfilled. 
Due to control redundancy, isolating faults affecting single actuators or effectors in overactuated systems can be a difficult task, as the same effects can be produced by faults occurring in different actuators or effectors: the family of filters needed to isolate the faults usually results to be larger compared to a control system framework with full-rank input matrix, and moreover there is an upper bound on the maximum number of simultaneously isolable faults. On the other hand, by constraining the inputs in prescribed configurations without altering system dynamics, control redundancy can be very helpful in separating the effects produced by multiple faults in order to identify which groups of effectors and actuators are losing effectiveness. \\
The paper is structured as follows. In Section 2 the basic setup of control allocation is introduced and the general structure of unknown input observers (UIO) is defined;  moreover some issues related to control reconfiguration are reported, such as control re-allocation in the presence of cost functions and input constraints. Section 3 is devoted to the presentation of the proposed method for designing families of detection/isolation filters based on UIO, namely {\it constrained output fault direction} (COFD) observers; an alternative approach based on a different class of unknown input observers has recently been proposed by the authors \cite{noiAut}. Finally, in Section 4,  the application of the theoretical results is extensively illustrated by the case study of an overactuated ship subject to thruster failures including common mode faults in thruster auxiliary systems. 
\section{Control allocation setup}
Let us consider the following linear system
$$
\left\{
\begin{array}{l}
\dot x(t)=A x(t)+B\tau(t)\\ \\
y(t)=C x(t)
\end{array}
\right.
$$
with
$$
\tau(t)=G u(t),
$$
where $x\in\mathbb{R}^n$, $\tau\in\mathbb{R}^k$, $y\in\mathbb{R}^p$, $u\in\mathbb{R}^m$, $m>k$ and all matrices except $A$ are assumed to be full-rank. The vector $x$ is the state, which is assumed to be not accessible for direct measurements, while $y$ is the measured output of the system. The vector $u(t)$ represents the redundant control input and $\tau(t)$ is the generalized control effect or virtual input. Without loss of generality, the desired control effect $\tau_c(t)$ is assumed to be given by a suitable known function depending on the system output:
\be\label{tau-y}
\tau_c(t)=f(y(t)).
\eeq
The above condition can also be generalized, assuming that the desired effect $\tau_c$ and the measured output $y(t)$ are related through a suitable dynamic law. 
A control allocation strategy is defined such that, whenever it is possible, the control $u$ satisfies
\be\label{unominal}
G u(t)=\tau_c(t).
\eeq
Although the above linear equation always admits (uncountable) exact solutions when $rk(G)=k$, there are possible constraints or bounds to be met and this may lead to the existence of approximate solutions only: 
\be\label{unominalcon}
\left\{\begin{array}{l}
u\in\mathbb{U}\\ 
Gu(t)=\hat\tau_c(t)\neq \tau_c(t)
\end{array}\right.
\eeq
We point out that the approximate effect $\hat\tau_c(t)$ may differ from the desired one $\tau_c(t)$ but it is a known quantity, as it can be computed exploiting the input constraints given by $\mathbb{U}$ which is an assigned set.\\ In the unconstrained case, a simple solution can be obtained using the right pseudo-inverse matrix \cite{Golub:83}:
\be\label{uncons}
u(t)=G^{-R}\tau_c(t),\quad G^{-R}:=G^T(GG^T)^{-1}.
\eeq
In the case of constrained control $u\in\mathbb{U}$, several methods for control allocation are available in the literature (\cite{Durham:93}, \cite{Shi:10}, \cite{BuEn:96}, \cite{Ha:04}, \cite{Za:09}, \cite{Bodson:02}, \cite{Frost:10}).
In this paper we consider the class of faults acting on effectors and actuators efficiency by changing their effectiveness: these can be modeled by a multiplicative term $\Delta(t)$:
$$
\tau(t)=G\Delta(t)u(t),\ \Delta(t)=diag[\delta_1(t),...,\delta_m(t)],
$$
for some unknown functions $\delta_i(t)$.
It follows that, whenever $\delta_i(t)\equiv1\ \forall i=1,...,m$, the controller operates with nominal conditions and hence
$$
\tau(t)=G u(t)=\hat\tau_c(t)
$$
On the other hand if one of the actuators is subject to a loss of effectiveness or complete failure, i.e. if $\delta_i(t)\neq1$ for some $i$, the designed control law will no longer be able to ensure the desired effect, this meaning that, in the case of fault presence, one may have
$$
\tau(t)\neq \hat\tau_c(t)
$$
with a consequent deterioration of system performances. Such problems can be avoided by accommodating the fault effects if a suitable control reconfiguration policy is considered. Defining a set of diagnosis signals, usually called residuals, one can detect and isolate the faults; then, performing the correct reconfiguration of the control input, one can track  the (approximate) desired effect $\hat\tau_c$ properly again. The approach presented in this paper is based on Unknown Input Observers UIO (see for instance \cite{CPZ:96}), whose general structure is the following:
$$
\left\{
\begin{array}{l}
\dot z(t)=Fz(t)+RB v(t) +Ky(t)\\ \\
\hat x(t)=z(t)+Hy(t)  
\end{array}
\right.
$$
where the matrices $F,R,K$ and $H$ are design parameters and $v(t)$, $y(t)$ are, respectively, a known reference input signal and the measured output of the system to be estimated through the observer; the signal $v(t)$ is usually set equal to the nominal and unperturbed input that is commanded to the system. It is worth to note that, in order to achieve a correct asymptotic state estimation, the matrix $F$ has to be Hurwitz.\\
As it will be shown in the following, unknown input observers are useful for the task of isolating faults in overactuated systems. Moreover, thanks to input redundancy, the control can be re-allocated in order to limit or avoid the use of faulty effectors/actuators once these have been isolated. On the other hand, this is not the only advantage of input redundancy in the considered framework: indeed, control allocation can be combined together with the fault isolation scheme in order to enlarge the family of identifiable faulty events.
%

The reconfiguration can be performed by different methods, depending on several factors such as actuator dynamics, bounds on energy consumption, limited control inputs rates or other control constraints.  It is worth to note that, due to the negative effects of faults, also the desired control effect $\tau_c(t)$ might be requested to change with respect to the original one in order to recover the deteriorated system performances, this corresponding to update the relation between $\tau_c(t)$ and the output signal $y(t)$ given by  (\ref{tau-y}).\\
In the simplest case of unconstrained inputs, the nominal control allocation law is given by (\ref{uncons}) and 
therefore, if the actuators $i_1,...,i_q$ are faulty and $q\leq m-k$, to get the desired effect $\tau_c(t)$ it is sufficient to re-allocate the control action setting 
$$
u_{i_1}=u_{i_2}=\cdots=u_{i_q}\equiv0
$$
and assigning the other components of $u(t)$, which are grouped for convenience in a vector $\tilde u\in\mathbb{R}^{m-q}$, according to 
\be\label{simple-realloc}
\tilde{u}(t)=\tilde G^{-R}\tau_c(t),
\eeq
where the matrix $\tilde G\in\mathbb{R}^{k\times (m-q)}$ is obtained from $G$ by neglecting the columns $i_1,...,i_q$.\\\\
On the other hand, since control reconfiguration can be regarded as a reduced-order control allocation problem in which some of the inputs are neglected, the use of the aforementioned techniques for handling input constraints can be straightforward extended. However, by turning off the input signals corresponding to faulty actuators, the redundancy of control inputs is reduced and the error between the desired control effect and the control effect provided by the approximate solution may increase after reconfiguration, as the class of admissible solutions to the allocation problem reduces.\\
As already mentioned, control allocation can be used actively also to make faults observable; in particular, by considering additional input constraints (see Section \ref{club}) which constrain control devices and surfaces to achieve common modes, one can isolate faults affecting selected groups of effectors. Whenever such constraints are not allowed to be imposed simultaneously in practice due to lack of control design freedom, an iterative control allocation scheme can be defined in order to switch periodically from one common mode to another after a prescribed time interval, until the fault isolation task is accomplished successfully.

\section{Fault detection and isolation}
The estimation error is defined as the difference between the true state $x(t)$ and the estimated state $\hat x(t)$:
$$
e(t)=x(t)-\hat x(t).
$$
Our aim is to design a family of unknown input observers $\{\mathcal{O}_h\}_{h=1}^s$, such that the information provided by the estimation errors allow us to detect and isolate faults. To address such target one can proceed as follows \cite{CPZ:96}.
\subsection{Unknown Input Observers}
The input $v(t)$ in the observer is set equal to the reference control effect $\hat\tau_c(t)$, that is
$$
v(t)=Gu(t)=\hat\tau_c(t),
$$
where $u(t)$ is the nominal (fault free) control (\ref{unominalcon}) (or (\ref{unominal}) in the unconstrained case). Exploiting the observer structure, the dynamics of the error is ruled by the following equation
$$
\begin{array}{ll}
\dot e(t)&=\dot x(t)-\dot{\hat x}(t)\\ 
&=[(I_{n\times n}-HC)A-KC+FHC]x(t)-F\hat x(t)\\ 
&+(I_{n\times n}-HC)BG\Delta(t)u(t)-RBG u(t).
\end{array}
$$
Setting $K=K_1+K_2$, if the following conditions are satisfied
\begin{align}
&R=I_{n\times n}-HC\label{uio1}\\ 
&F=RA-K_1C,\quad \sigma(F)\in\mathbb{C}^{-}\label{uio2}\\
&K_2=FH\label{uio3}
\end{align}
then the latter equation reduces to
$$
\dot e(t)=F e(t)+RBG(\Delta(t)-I_{m\times m})u(t),
$$
where $\sigma(\cdot)$ stands for the spectrum of a matrix and the set $\mathbb{C}^{-}$ in the left open complex half-plane. 
Let us denote by $W\in\mathbb{R}^{n\times m}$ the matrix $BG$, whose columns will be indicated with $W_1,...,W_m$, i.e.
\be\label{BG=W}
BG=W=[W_1\cdots W_m].
\eeq
It is worth to note that the matrix $BG(\Delta(t)-I)=W(\Delta(t)-I_{m\times m})$ appearing in the expression of $\dot{e}(t)$ has the following structure:
$$
W(\Delta(t)-I_{m\times m})=[(\delta_1(t)-1)W_1\ \cdots\ (\delta_m(t)-1)W_m],
$$and hence a fault in the $j^{th}$ effector may only affect the $j^{th}$ column of $W$.

\subsection{Constrained Output Fault Directions (COFD)}\label{prescrib}
This subsection is dedicated to the presentation of the proposed scheme for fault detection and isolation in overactuated systems; such method consists in constraining the residuals in prescribed subspaces of the output space (see for instance \cite{PaRi:94}, \cite{WhSp:87}).\\ Let us consider the canonical basis of $\mathbb{R}^p$, namely $\mathbf{e}_1,...,\mathbf{e}_p$; since by assumption the output matrix $C$ is full-rank, there exists $S\in\mathbb{R}^{n\times p}$ such that
$$
CS=I_{p\times p}=[\mathbf{e}_1\ \cdots\ \mathbf{e}_p].
$$
The general solution of such equation is given by
\be\label{SS0}
S=C^T(CC^T)^{-1}+[I_{n\times n}-C^T(CC^T)^{-1}C]S_*,
\eeq
where $S_*\in\mathbb{R}^{n\times p}$ is an arbitrary matrix. Denoting by $S_1,...,S_p$ the columns of the matrix $S$, the basic idea of the method is to design the observer parameters in order to guarantee that, if a fault occurs in the $i^{th}$ actuator, then the estimation error maintains the direction $S_i$ during the system evolution, this corresponding to a fixed direction $\mathbf{e}_i$ for the residual. It is worth to note that a first strong design constraint for the achievability of this condition is that directions $S_1,...,S_p$ need to correspond to eigenvectors of the observer matrix $F$. Moreover, due to the rank deficiency of $BG$, it is not possible in general to address a decoupled distribution of the faults effects over the columns of the matrix $S$ and we are required to deal with linear combinations of such characteristic directions.
We recall that the dynamics of the estimation error is given by the equation $$
\dot e(t)=F e(t)+RW(\Delta(t)-I_{m\times m})u(t).
$$
Since $rank(W)=k<m\leq p$, we can arbitrarily assign only $k$ columns of the matrix $RW$ through the design parameter $R$, as the remaining $m-k$ are consequently constrained; we need therefore to consider several independent observers to achieve a correct fault isolation. One can proceed as follows. 
\begin{deff}
We call a {$multi-index$} any vector $J$ of increasing natural numbers, i.e. $J=(j_1,...,j_\ell)$ with $ j_q\in\mathbb{N}\ \forall q=1,...,\ell$ and $1\leq j_1<j_2<\cdots< j_\ell\leq r,\  r\geq\ell$. The positive integers $\ell$ and $r$ are defined as, respectively, the length $L(J)$ and the domain $D(J)$ of multi-index $J$. 
\end{deff}
It is worth to note that the number $s$ of distinct multi-indices having length $\ell$ and domain $r$ is given by the  binomial coefficient 
\be\label{s-binom}
s=\binom{r}{\ell}=\frac{r!}{(r-\ell)!\ell!}.
\eeq
\begin{deff}\label{def2}
Given $W\in\mathbb{R}^{n\times m}$ with $rank(W)=k$, we call uniform sub-rank of $W$ the positive integer $k_0\leq k$ computed as
\be\label{kappa0}
\begin{array}{ll}
k_0:=\max&\!\!\!\!\!\{\ell\leq k: rank[W_{j_1}\ \cdots\ W_{j_\ell}]=\ell,\\
 &\!\!\forall\ J=(j_1,...,j_\ell): D(J)=m \}.
\end{array}
\eeq
\end{deff}
{\bf Notation} {\it Given a multi-index $J=(j_1,...,j_\ell)$ with $D(J)=m$, we denote by $W_J\in\mathbb{R}^{n\times\ell}$ the matrix composed by the columns of $W$ corresponding to the indices included in $J$, i.e.}
$$
W_J=[W_{j_1}\ \cdots\ W_{j_\ell}].
$$
Let $k_0\leq k$ the uniform sub-rank of the matrix $W=BG$. Then, if by tuning the matrix $R$ we prescribe the first $k_0$ columns of $RW$, for example imposing that they have to be equal to $S_1,S_2,...,S_{k_0}$, we get
$$
RW=[S_1\ \cdots\ S_{k_0}\ V_1\ \cdots\ V_{m-k_0}],
$$
where $V_j=\sum_{i=1}^{k_0}\alpha_{ij}S_i$ for some coefficients $\alpha_{ij}$. As a consequence we have
$$
CRW=[\mathbf{e}_1\ \cdots\ \mathbf{e}_{k_0}\  \omega_{1}\ \cdots\ \omega_{m-k_0}],
$$
with $\omega_j=\sum_{i=1}^{k_0}\alpha_{ij}\mathbf{e}_i.$ Setting $R^{(1)}=R$, we can iterate this construction by designing $R^{(h)}$ such that
\be\label{2met}
R^{(h)}[W_{j_1^h}\ \cdots\ W_{j_{k_0}^h}]=[S_{i_1^h}\ \cdots\ S_{i_{k_0}^h}]
\eeq
as multi-indices $J_h=(j_{1}^h,...,j_{k_0}^h)$ and $I_h=(i_{1}^h,...,i_{k_0}^h)$ vary; at the end of the construction we obtain a family of matrices $\{R^{(h)}\}_{h=1}^s$, with
$$
s=\binom{m}{k_0}\binom{p}{k_0}=\frac{m!p!}{(m-k_0)!(p-k_0)!(k_0!)^2}.
$$
We point out that the information provided by the residuals associated to such family of matrices is redundant. In order to reduce the computational burden and avoid overlapping of information, we need to investigate what is the minimum number of matrices $R^{(h)}$ required for a proper fault isolation. We see from the explicit construction of $R^{(1)}$ that, if $k_0>1$, in this case faults affecting the first $k_0$ actuators lead to residual signals directed as $\mathbf{e}_1,...,\mathbf{e}_{k_0}$ respectively, while faults affecting the other actuators lead to residual signal obtained as linear combination of two or more vectors $\mathbf{e}_j$, $j=1,...,k_0$. On the other hand, if more than one fault occurs we get a residual signal defined by a linear combination of vectors $\mathbf{e}_j$ as well, this meaning that with the information provided by this unique observer we are not able to distinguish multiple faults from individual faults affecting one of the last $m-k_0$ actuators. For this reason, the number of observers to be considered has to be sufficient to decouple effects of single and multiple faults. To this purpose we note that, since the observers are designed independently, the particular choice of vectors 
$S_j$ among $\{S_1\ \cdots\ S_p\}$ in (\ref{2met}) does not influence the fault isolation procedure: nevertheless, this freedom of choice may result to be helpful for the design of the Hurwitz matrix $F$. We claim that the maximum number of isolable faults is $k_0-1$ with a required number of observers equal to
$$
\bar{s}=\binom{m}{k_0}.
$$
This can be verified observing that, in order to isolate the faults, at least one residual $r^{(h)}(t)$ needs to have null projection along one of the basis vectors $\mathbf{e}_j$: as a consequence, the number of isolable faults has to be less than the maximum admissible number of independent components of each residual $r^{(h)}(t)$, that is $k_0$. The integer $\bar{s}$ can be obtained simply computing the number of all structurally distinct vectors of $m$ elements with $k_0$ assigned entries.\\
For sake of simplicity we fix the vectors $S_j$ in (\ref{2met}), assuming that the right-hand side is equal to $\hat{S}:=[S_1\ \cdots\ S_{k_0}]$ for any $h$.
Once the properties of the matrices $R^{(h)}$ are defined, one have to deal with the stabilization of the matrix $F^{(h)}$ together with the fulfillment of the rank condition
\be\label{rank1}
rank[S_j\ F^{(h)}S_j\ \cdots\ (F^{(h)})^{n-1}S_j]=1\ \forall\ j=1,...,k_0,
\eeq
which corresponds to the requirement for $S_j$ to be an eigenvector of the matrix $F^{(h)}$, this last condition being fundamental for ensuring a constant direction of the output.
We point out that for any $h=1,...,\bar{s}$, the solution of (\ref{2met}) is given by
\be\label{acca}
\left\{
\begin{array}{ll}
R^{(h)}&=(I_{n\times n}-H^{(h)}C)\\
H^{(h)}&=(W_{J_h}-\hat S)(CW_{J_h})^{-L}\\
&+H_*^{(h)}(C-(CW_{J_h})(CW_{J_h})^{-L}C)
\end{array}\right.
\eeq
where $(\cdot)^{-L}$ stands for the left pseudo-inverse and $H^{(h)}_*\in\mathbb{R}^{n\times p}$ is an arbitrary matrix. Let $M^{(h)}\in\mathbb{R}^{k_0\times k_0}$ be a diagonal and negative definite matrix. Recalling that $F^{(h)}=R^{(h)}A-K_1^{(h)}C$, condition (\ref{rank1}) can be achieved by solving
\be\label{gain}
K_1^{(h)}(C\hat S)=R^{(h)}A\hat S-\hat S M^{(h)},
\eeq
which gives
\be\label{gain2}
\begin{array}{ll}
K^{(h)}_1&=(R^{(h)}A\hat S-\hat S M^{(h)})(C\hat S)^{-L}\\ 
&+K_*^{(h)}(I_{p\times p}-(C\hat S)(C\hat S)^{-L}),
\end{array}
\eeq
where $K_*^{(h)}\in\mathbb{R}^{n\times p}$ is an arbitrary matrix. Denoting by $e^{(h)}(t)$ the estimation error associated to the observer $\mathcal{O}_h$, we define the (vectorial) residual signals:
$$
r^{(h)}(t)=Ce^{(h)}(t).
$$
 We are now ready to state the main result of this section.
\begin{theo}\label{theomain2}
Let $S, R^{(h)}, H^{(h)}, K^{(h)}_1$ be assigned by (\ref{SS0}), (\ref{acca}) and (\ref{gain2}).
Let us assume that, for any multi-index $J_h=(j_1^h,...,j^h_{k_0})$, $h=1,...,\bar s$, the following conditions hold true
\begin{enumerate}
\item $rank(W_{J_h})=rank(CW_{J_h})=k_0$;
\item the matrices $S_*, H^{(h)}_*$ and $K_*^{(h)}$ can be found such that the matrix $F^{(h)}=R^{(h)}A-K_1^{(h)}C$ is Hurwitz.
\end{enumerate}
Then $(F^{(h)}, H^{(h)}, R^{(h)}, K^{(h)})$ is an unknown input observer and it satisfies condition (\ref{2met}) for any $h=1,...,\bar s$, hence the residual signals are able to detect and isolate up to $k_0-1$ faults affecting simultaneously the system actuators. 
\end{theo}

We can represent residuals as ordered sums of the basis vectors and their combinations; we will employ the $\oplus$ to indicate a logic sum of the vectors $v,w$ depending on the order, i.e. $v\oplus w\neq w\oplus v.$ For example we obtain the following logic representation of the first residual:
$$
r^{(1)}=\mathbf{e}_1\oplus\mathbf{e}_2\oplus\cdots\oplus\mathbf{e}_{k_0}\oplus\omega_1^1\oplus\cdots\oplus\omega_{m-k_0}^1
$$
where $\omega_j^1$ is an arbitrary combination of the vectors $\mathbf{e}_1,...,\mathbf{e}_{k_0}$; in a similar way, the residual associated to the multi-index $J_h=(j_1^h,...,j_{k_0}^h)$ can be represented by
$$
r^{(h)}=\omega_{1}^h\oplus\cdots\oplus\omega_{j_1-1}^h\oplus\mathbf{e}_1\oplus\omega_{j_1}^h\oplus\cdots\oplus\mathbf{e}_{k_0}\oplus\cdots\oplus\omega_{m-k_0}^h
$$
\begin{rem}
It is worth to note that condition (\ref{rank1}) is sufficient for prescribing a fixed direction to the residual signals only if  the observer initialization error is zero or if the estimation error $e(t)$ has reached the steady-state at the moment of fault occurrence.
\end{rem}
\subsection{Cluster residuals}\label{club}
We consider here an extended framework, in which the actuators are grouped into $q$ clusters:
$$
\begin{array}{l}
\mathcal{A}_1=\{u_1,...,u_{i_1}\},\ \mathcal{A}_2=\{u_{i_1+1},...,u_{i_2}\},\\
\ \dots\dots\ \mathcal{A}_q=\{u_{i_{q-1}+1},...,u_m\},
\end{array}
$$
where $1\leq i_1<i_2<\cdots i_{q-1}<m$. The introduction of this new model is motivated by the need of isolate common mode faults which may affect simultaneously actuators or effectors sharing the same auxiliary systems. For sake of clarity only the case of non-overlapping clusters is presented; on the other hand, the proposed methods can be readily modified in order to be used also in the case of possibly overlapping clusters. \\
Let us denote by $\alpha_h$ the cardinality of the cluster $\mathcal{A}_h$, i.e. $\alpha_h=i_h-i_{h-1}$ and hence $\sum_{i=1}^q\alpha_i=m$.
The faults are supposed to act uniformly on the whole cluster $\mathcal{A}_i$, so that they can be modeled as the block-diagonal multiplicative matrix:
\be\label{deltaclust}
\Delta(t)=diag(d_1(t)I_{\alpha_1\times \alpha_1},d_2(t)I_{\alpha_2\times\alpha_2},\cdots,d_q(t)I_{\alpha_q\times\alpha_q})
\eeq
Let $k_0$ be the uniform sub-rank of $W$ and let us suppose that $k_0\geq3$ and
$$
\max_{i=1,...,q}\alpha_i\leq k_0-1.
$$
As a consequence if a fault is present in a single cluster it can be detected and isolated; in the same way, if there exist two (or more clusters) with $\alpha_i+\alpha_j\leq k_0-1$, a multiple fault can be isolated as well. On the other hand, if for some pair of indices $i,j$, one has $\alpha_i+\alpha_j\geq k_0$, a multiple fault on the associated effectors will lead to completely saturated residuals, this meaning that no fault isolation can be achieved at the present step. However, by introducing additional constraints in the control allocation scheme, one can design a new set of observers to be used to identify faulty clusters of actuators/effectors. For sake of simplicity let us consider first the following case:
$$
\alpha_1+\alpha_2\geq k_0;
$$
in addition, let us assume that
\be\label{m-k}
2+m-k\geq\alpha_1+\alpha_2.
\eeq
We select two finite sequences of real numbers $\{\zeta_1^{(1)},\cdots,\zeta_{\alpha_1-1}^{(1)}\}$ and $\{\zeta_{1}^{(2)},\cdots,\zeta^{(2)}_{\alpha_2-1}\}$; using the control redundancy ensured by (\ref{m-k}), we are allowed to impose the constraints
\be\label{redund}
\begin{array}{l}
\displaystyle\frac{u_j(t)}{u_{i_1}(t)}=\zeta^{(1)}_{j},\quad j=1,...,\alpha_1-1,\\ \\
\displaystyle\frac{u_{i_1+j}(t)}{u_{i_2}(t)}=\zeta^{(2)}_{j},\quad j=1,...,\alpha_2-1
\end{array}
\eeq
together with the equality
$$
\tau_c(t)=Gu(t).
$$
Due to (\ref{redund}), the overall input signal associated to the first cluster $\mathcal{A}_1$ turns out to be
$$
\sum_{i=j}^{i_1}W_ju_j(t)=u_{i_1}(t)\left(W_{i_1}+\sum_{j=1}^{\alpha_1-1}{\zeta^{(1)}_j}W_j\right),
$$
while the overall signal corresponding to $\mathcal{A}_2$ is
$$
\sum_{j=i_1+1}^{i_2}W_ju_j(t)=u_{i_2}(t)\left(W_{i_2}+\sum_{j=1}^{\alpha_2-1}{\zeta^{(2)}_j}W_{j+i_1}\right).
$$
Let us denote by $W^{\{1\}}$ and $W^{\{2\}}$ the constant vectors
$$
W_{\{1\}}=W_{i_1}+\sum_{j=1}^{\alpha_1-1}{\zeta^{(1)}_j}W_j,\ W_{\{2\}}=W_{i_2}+\sum_{j=1}^{\alpha_2-1}{\zeta^{(2)}_j}W_{j+i_1};
$$
since, without loss of generality, the coefficients $\zeta_j^{(\star)}$, $\star=i_1,i_2$, can be chosen such that the latter vectors are independent, by definition one has
\be\label{rango3}
rank[W_{\{1\}}\ W_{\{2\}}\ W_{i_2+1}\ \cdots\ W_{m}]\geq 3.
\eeq
If the cluster $\mathcal{A}_1$ undergoes a fault, the dynamics of the estimation error turns out to be
$$
\dot e(t)=F e(t)+ RW_{\{1\}}(\delta_1(t)-1)u_{i_1}(t)
$$
and a similar condition holds for faults in $\mathcal{A}_2$. Now, recalling that $rank[W_{\{1\}}\ W_{\{2\}}]= 2$ and designing the observer matrix $R$ such that
$$
R[ W_{\{1\}}\ W_{\{2\}}\ W_i]=[\mathbf{e_1}\ \mathbf{e}_2\ \mathbf{e}_3],
$$
for some $i$ with $u_i\notin \mathcal{A}_1\cup \mathcal{A}_2$, one obtains a residual signal with prescribed output directions associated to the first two clusters of actuators. We will refer to such signal as a {\it cluster residual}. The above construction can be readily extended to the case of faults involving more than two effector clusters, if there is enough redundancy to use control allocation. To this purpose, let us set
\be\label{wil}
W_{\{h\}}=W_{i_h}+\sum_{j=1}^{\alpha_h-1}\zeta^{(i_h)}_jW_{i_{h-1}+1},\quad h=1,...,q
\eeq
where the coefficients $\zeta^{(i_\ell)}_j\in\mathbb{R}$ have to be defined. Such overall input vectors can be organized into a reduced-order input matrix
$$
W^\star=[W_{\{1\}}\ W_{\{2\}}\ \cdots\ W_{\{q\}}];
$$
applying the observer design scheme proposed in Section III.A with $W$ replaced by $W^*$, one can obtain a family of cluster residuals associated to the effector groups $\mathcal{A}_1,...,\mathcal{A}_q$.
\section{Case study: overactuated marine vessel}
\begin{figure}[h!]
\centering
{\includegraphics[height=0.25\textheight]{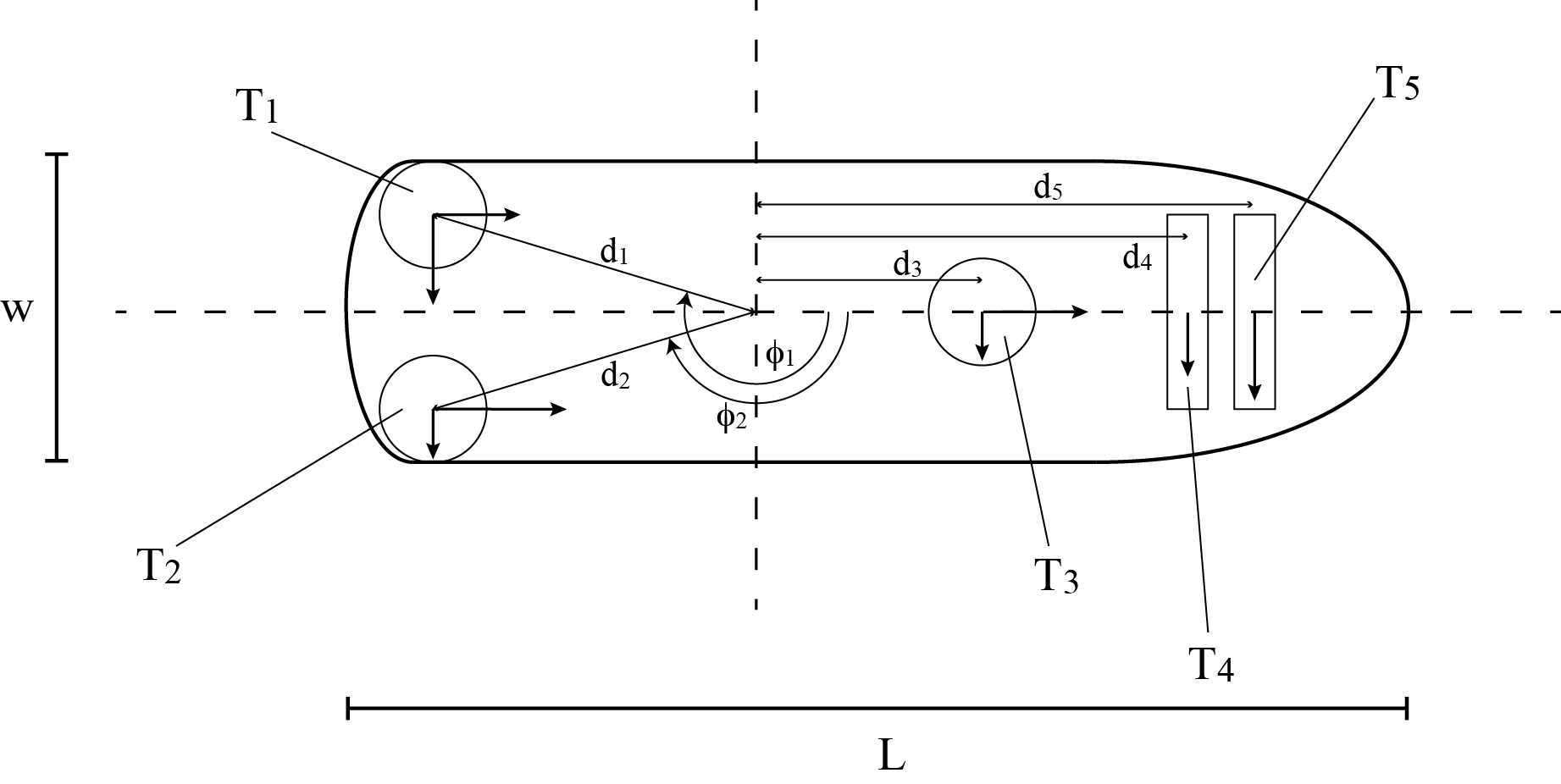}}
\caption{\label{fig1}Five thrusters ship model}
\end{figure}

This section is focused on illustrating the application of theoretical results to the case of an overactuated marine vessel. We consider the following ship model 
$$
\left\{
\begin{array}{l}
\dot\eta =  P(\eta)\nu\\ \\
M\dot\nu=-V(\eta,\nu)+\tau+P^T(\eta)b(t)
\end{array}
\right.
$$
where $M$ is the inertia matrix, $\eta=[x_G,y_G,\psi]^T$ is the ship position coordinates in the earth-fixed reference frame and $\nu=[\nu_x,\nu_y,\dot\psi]$ contains surge, sway and yaw angular velocities with respect to the body-fixed reference frame, which is identified with the vessel center of mass; the vector $V(\eta,\nu)$ includes Coriolis and damping terms, the actual thrust force in surge, sway and the yaw moment are given by $\tau=[\tau_x,\tau_y,m_{\dot\psi}]$ and $P(\eta)=P(\psi)$ is a standard rotation matrix
$$
P(\psi)=\left[
\begin{array}{rrr}
\cos\psi&-\sin\psi&0\\
\sin\psi&\cos\psi&0\\
0&0&1
\end{array}
\right]
$$
The perturbation term $b(t)$, which is assumed to be bounded by a known constant, is used to model disturbances affecting the system, such as slowly-varying forces and moments caused by wave loads, ocean currents or winds (see \cite{Fossen:94} for further details).
Following \cite{BeFo:97}, both $\eta$ and $\nu$ are assumed to be measured (possibly through a state estimation), the constant matrix $M$ is known and $V(\eta,\nu)$ is a known function in the state variables $\eta,\nu$. We consider a ship equipped with 3 azimuth thrusters (rotatable) $T_1, T_2, T_3$ and 2 transverse tunnel thrusters (fixed orientation) $T_4, T_5$. A sketch of the vessel model is depicted in Figure \ref{fig1}.
\\
The actual thrust force is related to the control input through the linear equation
$$
\tau=Gu(t),
$$
with
$$
G=\left[
\begin{array}{rrrrrrrr}
1&0&1&0&1&0&0&0\\
0&1&0&1&0&1&1&1\\
\gamma_1&\gamma_2&\gamma_3&\gamma_4&\gamma_5&\gamma_6&\chi_7&\chi_8
\end{array}
\right],
$$
where the moment arms $\gamma_j$ are associated to azimuth thrusters and the moment arms $\chi_j$ to tunnel thrusters instead. In particular such arms can be computed as
$$
\begin{array}{rll}
\gamma_{2s-1}&=d_s\sin{\phi_s}& s=1,2,3\\
\gamma_{2s}&=d_s\cos{\phi_s}& s=1,2,3\\
\chi_{s}&=d_s\cos\phi_s& s=4,5
\end{array}
$$
where $d_s$ are the distances of thrusters and $\phi_s$ are the angles from the rotation point.\\
Assuming that the vessel rotation is negligible with respect to translation motion, setting $X=[\eta,\nu]^T$, for sake of simplicity and without loss of generality, the above nonlinear model can be linearized as follows
\be\label{linearmod}
\dot X(t)=AX(t)+B\tau(t)+E(t)
\eeq
with 
$$
A=\left[
\begin{array}{rr}
0_{3\times 3}& P(\bar\psi)\\ 
0_{3\times 3}&-M^{-1}D
\end{array}
\right],\quad\quad B=\left[
\begin{array}{r}
0_{3\times 3}\\
M^{-1}
\end{array}
\right],
$$
where $\bar\psi$ is a constant angle associated to some reference heading direction, $D=D(\bar\nu)$ is a constant damping matrix depending on a nominal reference velocity $\bar\nu$ and
$$
E(t)=\left[
\begin{array}{r}
0_{3\times 1}\\
 M^{-1}P^T(\bar\psi)b(t)
\end{array}
\right].
$$
Since by assumption all state variables are measured, without loss of generality the output matrix $C$ is supposed to be equal to the identity matrix
$$
C=I_{6\times 6}.
$$
Assuming that the marine vessel has a mass $\mu=6\cdot10^{6}Kg$, with lenght $L=76 m$ and width $w= 16 m$, the following parameters are obtained \cite{Fossen:94}:
$$
\begin{array}{l}
M=10^9\left[
\begin{array}{rrr}
0.0068 & 0. & 0.\\
0. & 0.0113 & -0.0340\\
0. & -0.0340& 4.4524
\end{array}
\right],\\ \\ D=10^8\left[
\begin{array}{rrr}
0.0008 & 0. &0.\\
0 & 0.0025 & -0.0203\\
0 & -0.0340 & 3.8481
\end{array}
\right].
\end{array}
$$
Without loss of generality $\bar\psi=0$, that is $P(\bar\psi)=I_{3\times 3}$. Moreover, taking $d_1=d_2=20m$, $d_3=18.5m$, $d_4=30m$, $d_5=35m$ and $\phi_1=\pi+0.3,\ \phi_2=\pi-0.3$, the matrix $G$ is given by
$$
G\simeq\!\!\left[
\begin{array}{rrrrrrrr}
1&0&1&0&1&0&0&0\\
0&1&0&1&0&1&1&1\\
-5.91 & -19.1 & 5.91 & -19.1 & 0& 18.5 & 30 & 35
\end{array}
\right]\!\!.
$$
We suppose the faults to occur in effectors (thrusters) rather than in single actuators: this corresponds to consider a fault matrix $\Delta(t)=diag(\delta_1(t),\delta_1(t),\delta_2(t),\delta_2(t),\delta_3(t),\delta_3(t),\delta_4(t),\delta_5(t))$ with coupled entries representing failures in the thrusters $T_j$, $j=1,...,5$. We set
$$
W=BG=[W_1\ W_2\ \vline\ W_3\ W_4\ \vline\ W_5\ W_6\ \vline\ W_7\ \vline\ W_8 ],
$$
where the vertical rules have been added to easily individuate the actuators corresponding to each thruster.\\
We can design a family of $s=4$ COFD unknown input observers $\{\mathcal{O}_h\}$ to isolate faults affecting singularly each thruster. In particular, following the steps of Theorem \ref{theomain2}, we select the observers matrices $R^{(h)}$ in order to have
$$
\begin{array}{ll}
R^{(1)}[W_1\ W_2\ W_3]=[\mathbf{e}_1\ \mathbf{e}_2\ \mathbf{e}_3]\\
R^{(2)}[W_3\ W_4\ W_1]=[\mathbf{e}_1\ \mathbf{e}_2\ \mathbf{e}_3]\\
R^{(3)}[W_5\ W_6\ W_1]=[\mathbf{e}_1\ \mathbf{e}_2\ \mathbf{e}_3]\\
R^{(4)}[W_7\ W_8\ W_1]=[\mathbf{e}_1\ \mathbf{e}_2\ \mathbf{e}_3].
\end{array}
$$
The observer gains $K_1^{(h)}$ can be chosen in order to assign the eigenvalues of the observer matrices $F^{(h)}$; on the other hand, since the whole state $X$ is measurable, i.e. $C=I_{6\times 6}$, the matrices $F^{(h)}$ turn out to be diagonal. For sake of simplicity we assume $F^{(h)}=F$ for any $h=1,...,5$, with $F=diag(-1,-1,-2,-5,-6,-7)$.\\
Different single or multiple thruster fault events have been simulated, assuming the initial conditions $\eta_0=\eta(0)=(1,1,0)$ and $\nu_0=\nu(0)=(2.2,1.9,0)$. The disturbance term is supposed to be given by the sum of two contributions: a constant term with random but fixed input direction representing an irrotational ocean current and an oscillating term with varying input direction representing waves;  the overall disturbance effect $b(t)$ is assumed to be bounded by the known constant $\epsilon=5\cdot10^6$.
The ship is supposed to be equipped with an $xy$-joystick control device together with a heading autopilot; the nominal operating conditions of the vessel are defined by a constant translational speed regime, this corresponding to the commanded control effect
$$
\tau_c(t)=D\nu(t)+[0\ 0\ a_\psi]^T,
$$
where $a_\psi$ is a PID controller for the yaw angle.\\
\subsubsection*{Single faults} 
\begin{figure}
\centering
{\includegraphics[height=0.25\textheight]{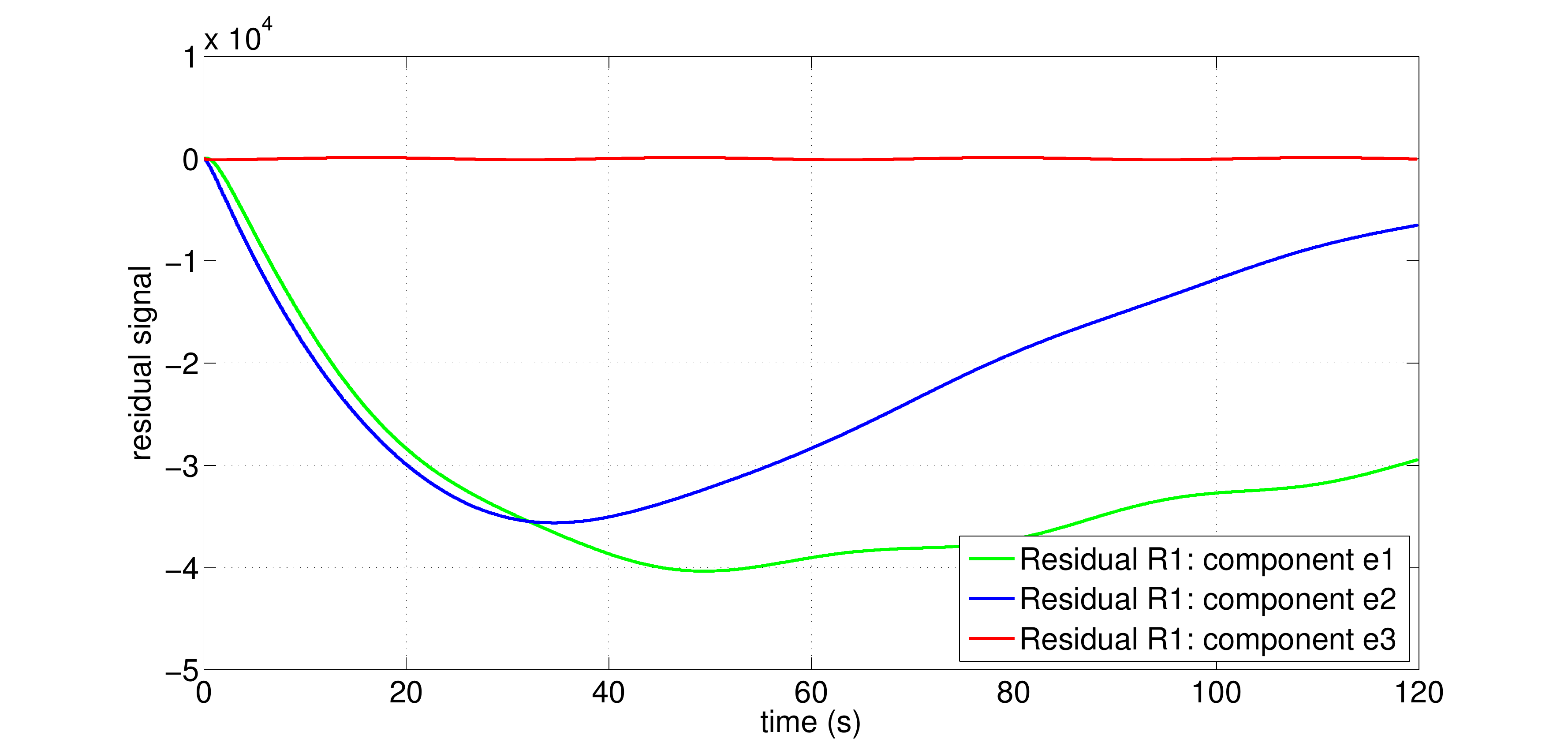}}
\caption{\label{fig2}Five thrusters ship model}
\end{figure}

\begin{figure}
\centering
{\includegraphics[height=0.25\textheight]{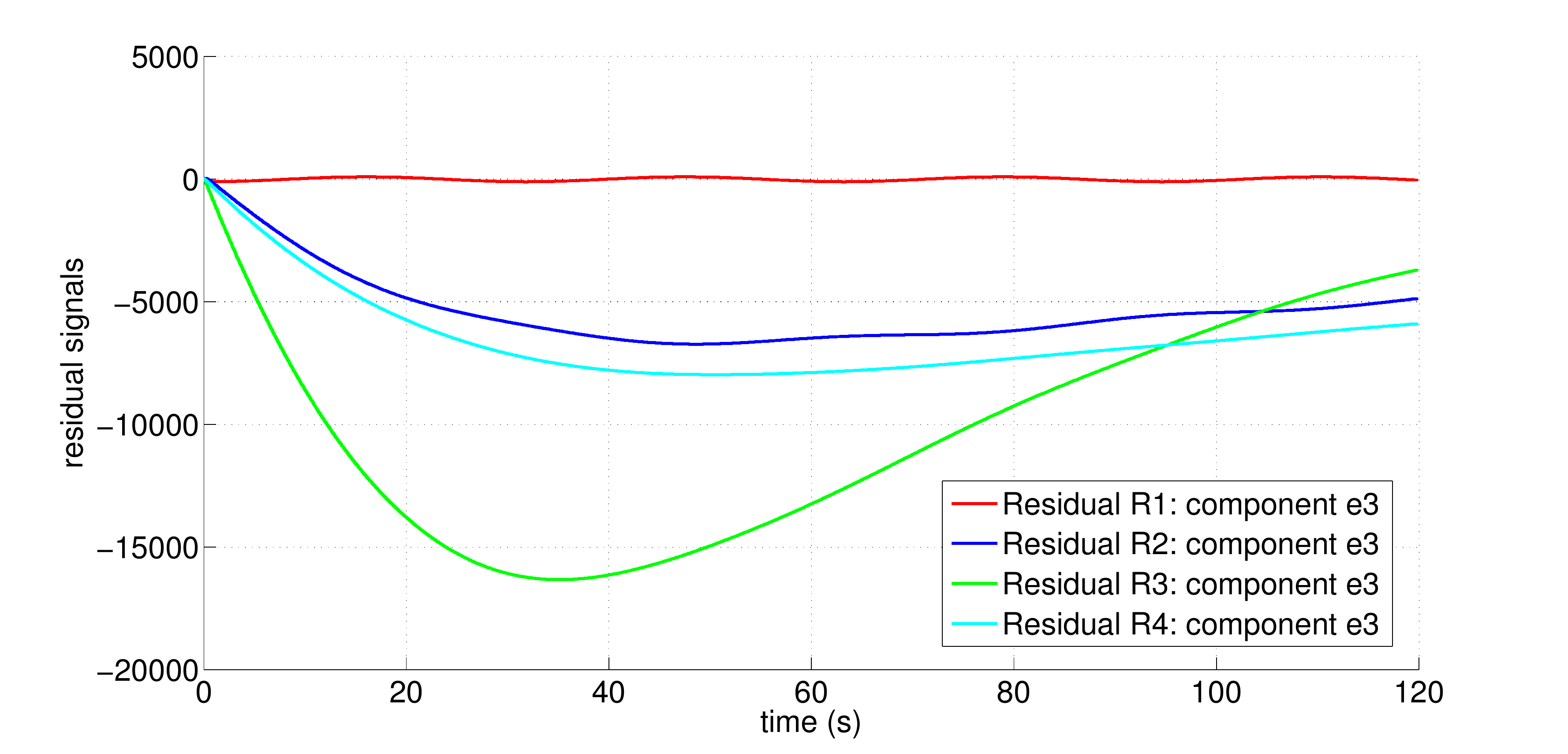}}
\caption{\label{fig3}Five thrusters ship model}
\end{figure}

We first suppose the azimuth thruster $T_1$ to be affected by a fault that gradually fades and its effect is $\delta_1(t)=e^{-0.03t}$. 
The behavior of the residual $r^{(1)}(t)$ is depicted in Figure \ref{fig2}: while the projections of $r^{(1)}$ along the directions $\mathbf{e}_1$ and $\mathbf{e}_2$ are significant, the projection along the direction $\mathbf{e}_3$ is negligible; on the other hand, as shown in Figure \ref{fig3}, the projections along the direction $\mathbf{e}_3$ of the other residuals is not negligible. 
As a consequence the faulty thruster $T_1$ can be identified and the control reconfiguration policy can be applied. \\
\\
\subsubsection*{Multiple faults}
\begin{figure}
\centering
{\includegraphics[height=0.25\textheight]{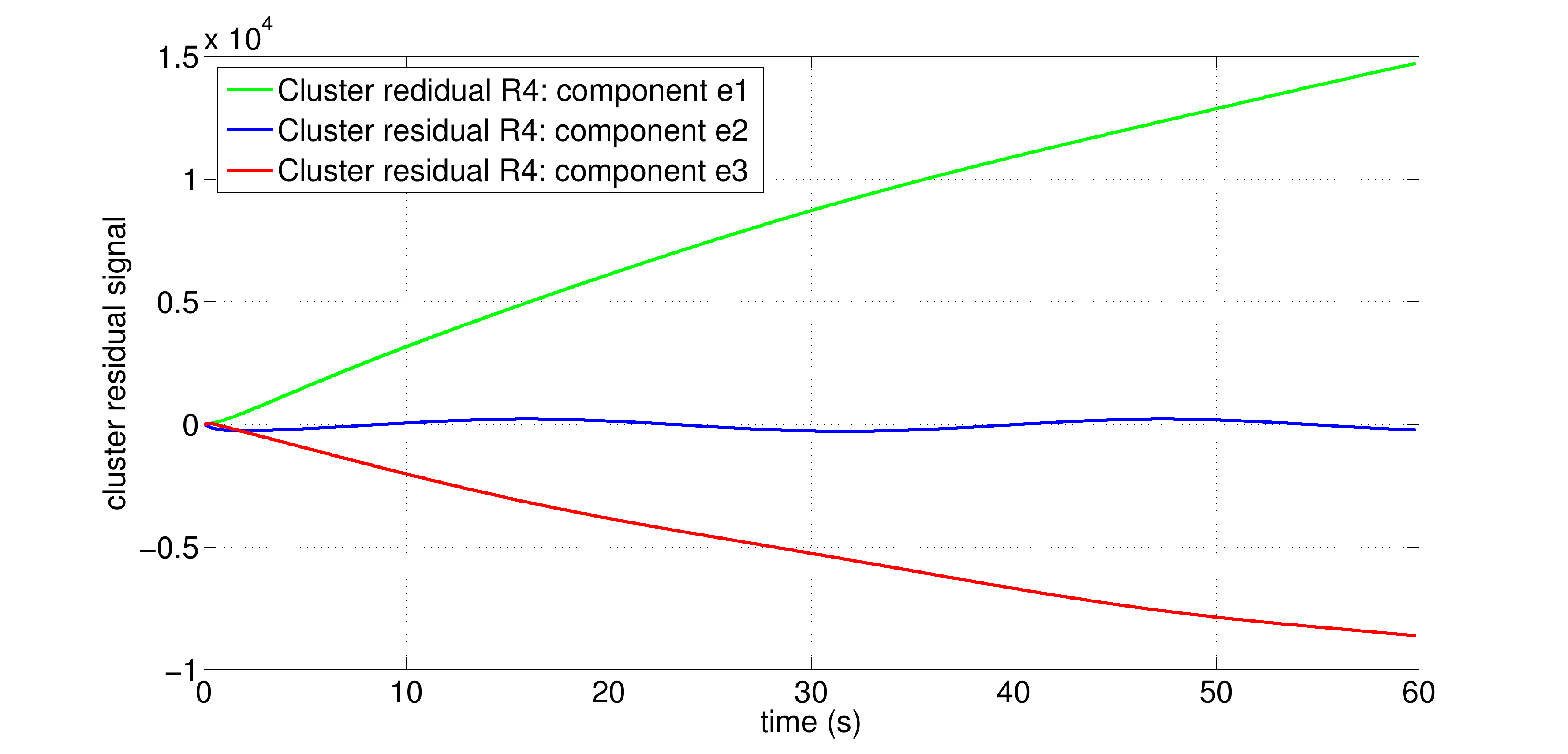}}
\caption{\label{fig8}Faults affecting simultaneously the thrusters $T_2$ and $T_5$: evolution of the cluster residual $R^{(4)}$}
\end{figure}

The most relevant cases of faults affecting simultaneously two thrusters due to common auxiliaries or power supply failures are reported in the following table:

\footnotesize{\begin{center}
    \begin{tabular}{ | c | c | c | c | c | c |}
    \hline
     & $T_1$ & $T_2$ & $T_3$ & $T_4$ & $T_5$ \\ \hline
     $T_1$ & & & $\bigcirc$ &$\bigcirc$ &\\ \hline
     $T_2$ & & & $\bigcirc$& &$\bigcirc$\\ \hline
     $T_3$ & $\bigcirc$&$\bigcirc$ & &$\bigcirc$ &$\bigcirc$\\ \hline
      
    \end{tabular}\\ \vspace{0.3cm}
 \scriptsize{Table 1 - The circles indicate that the case of a fault\\ occurring in the corresponding pair of thrusters is relevant}   
\end{center}}
\normalsize
In particular, it is assumed that $T_1$ and $T_4$ share the same auxiliaries, as well as $T_2$ and $T_5$ do; the thruster $T_3$ is supposed to be equipped with a switching device that enables it to be connected arbitrarily to  one sub-group or to the other, depending on the operating conditions of the system.\\ Using actively the control allocation, a new family of observers can be designed to obtain cluster residuals (as showed in Section \ref{club}); for sake of simplicity we assume the evolution time to be re-initialized at the present step. We choose the coefficients $\zeta^{(1)}=2.27, \zeta^{(2)}=3.41$ and $\zeta^{(3)}=1.38$ according to the initial values of the control inputs, i.e. $\zeta^{(i)}=u_{2i}(0)/u_{2i-1}(0)$; we impose the following conditions on the control input
\be\label{glizeta}
\zeta^{(1)}u_1(t)=u_2(t),\ \zeta^{(2)}u_3(t)=u_4(t),\ \zeta^{(3)}u_5(t)=u_6(t)
\eeq
and we set
\be\label{gliW}
\begin{array}{l}
W_{\{1\}}=\zeta^{(1)}W_1+W_2,\  W_{\{2\}}=\zeta^{(2)}W_3+W_4,\\
 W_{\{3\}}=\zeta^{(3)}W_5+W_6,\ W_{\{4\}}=W_7,\ W_{\{5\}}=W_8.
\end{array}
\eeq
By transformation (\ref{glizeta})-(\ref{gliW}), each thruster $T_i$ is associated to its overall input vector $W_{\{i\}}$, $i=1,2,3,4,5$, in a scheme that corresponds to fixed directions of azimuth thrusters. A family of COFD unknown input observers $\{\mathcal{O}_h\}_{h=1}^{4}$ is then defined such that the corresponding matrices $R^{(h)}$ satisfy:
\be\label{rules}
\begin{array}{l}
R^{(1)}[W_{\{1\}}\ W_{\{2\}}\ W_{\{3\}}]=[\mathbf{e}_1\ \mathbf{e}_2\ \mathbf{e}_3 ]\\
R^{(2)}[W_{\{1\}}\ W_{\{4\}}\ W_{\{5\}}]=[\mathbf{e}_1\ \mathbf{e}_2\ \mathbf{e}_3 ]\\
R^{(3)}[W_{\{2\}}\ W_{\{3\}}\ W_{\{4\}}]=[\mathbf{e}_1\ \mathbf{e}_2\ \mathbf{e}_3 ]\\
R^{(4)}[W_{\{2\}}\ W_{\{3\}}\ W_{\{5\}}]=[\mathbf{e}_1\ \mathbf{e}_2\ \mathbf{e}_3 ]
\end{array}
\eeq
Figure \ref{fig8} illustrates the case of faults $\delta_2(t)=e^{-0.02t }$ and $\delta_5(t)=e^{-0.01t}$ affecting simultaneously the thrusters $T_2, T_5$.  The cluster residual $r^{(4)}$ has a negligible component along the direction $\mathbf{e}_2$ while the other components are significant; moreover it can be verified that all other cluster residuals have significant components along the directions $\mathbf{e}_i$ for $i=1,2,3$. Based on the design rules (\ref{rules}), such information is sufficient to allow a correct identification of the pair of faulty thrusters $T_2$ and $T_5$.
\vspace{0.2cm}
\subsubsection*{Control reconfiguration}
\begin{figure}
\centering
{\includegraphics[height=0.25\textheight]{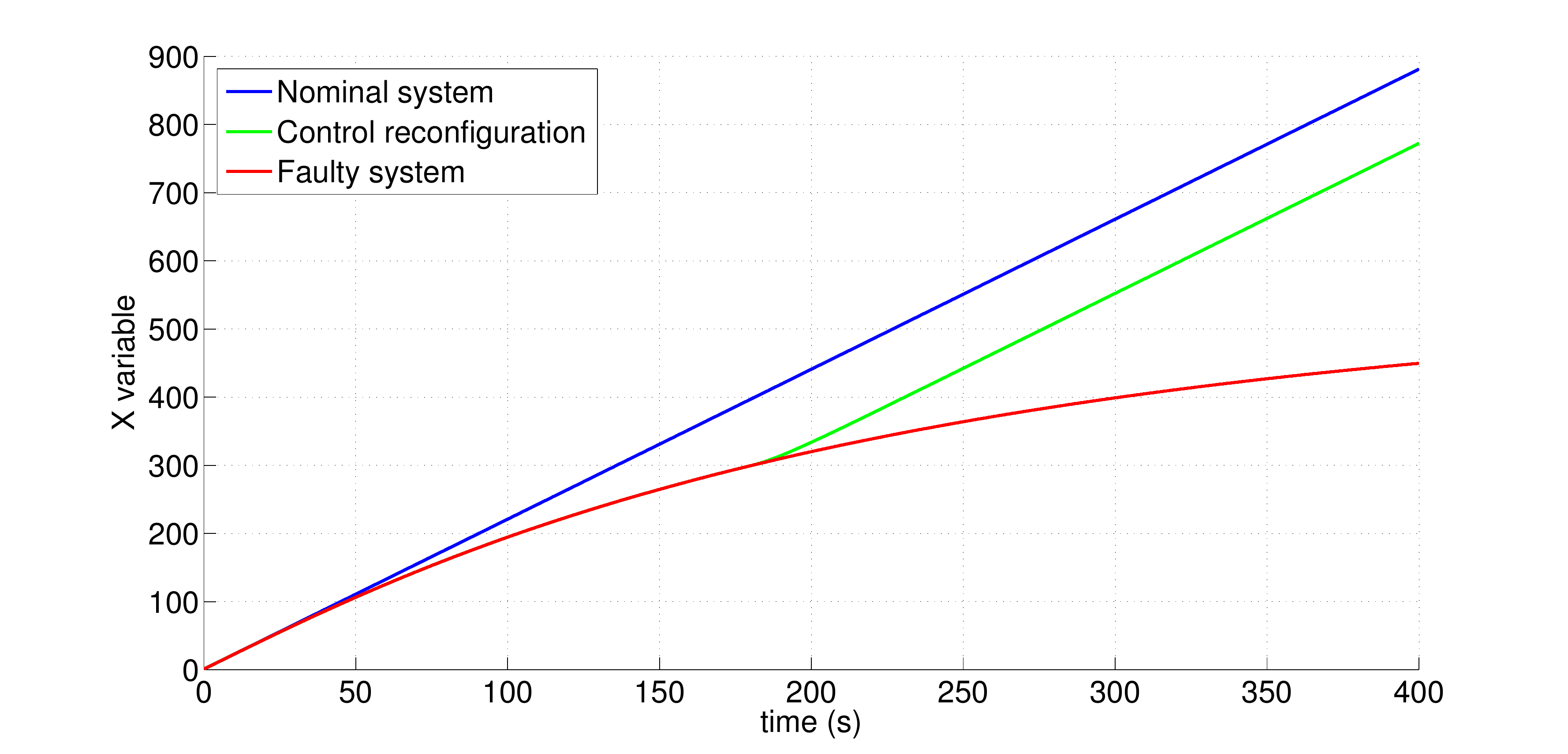}}
\caption{\label{fig5}Fault affecting the thruster $T_1$: evolution of the state variable $x_G(t)$}
\end{figure}

\begin{figure}
\centering
{\includegraphics[height=0.25\textheight]{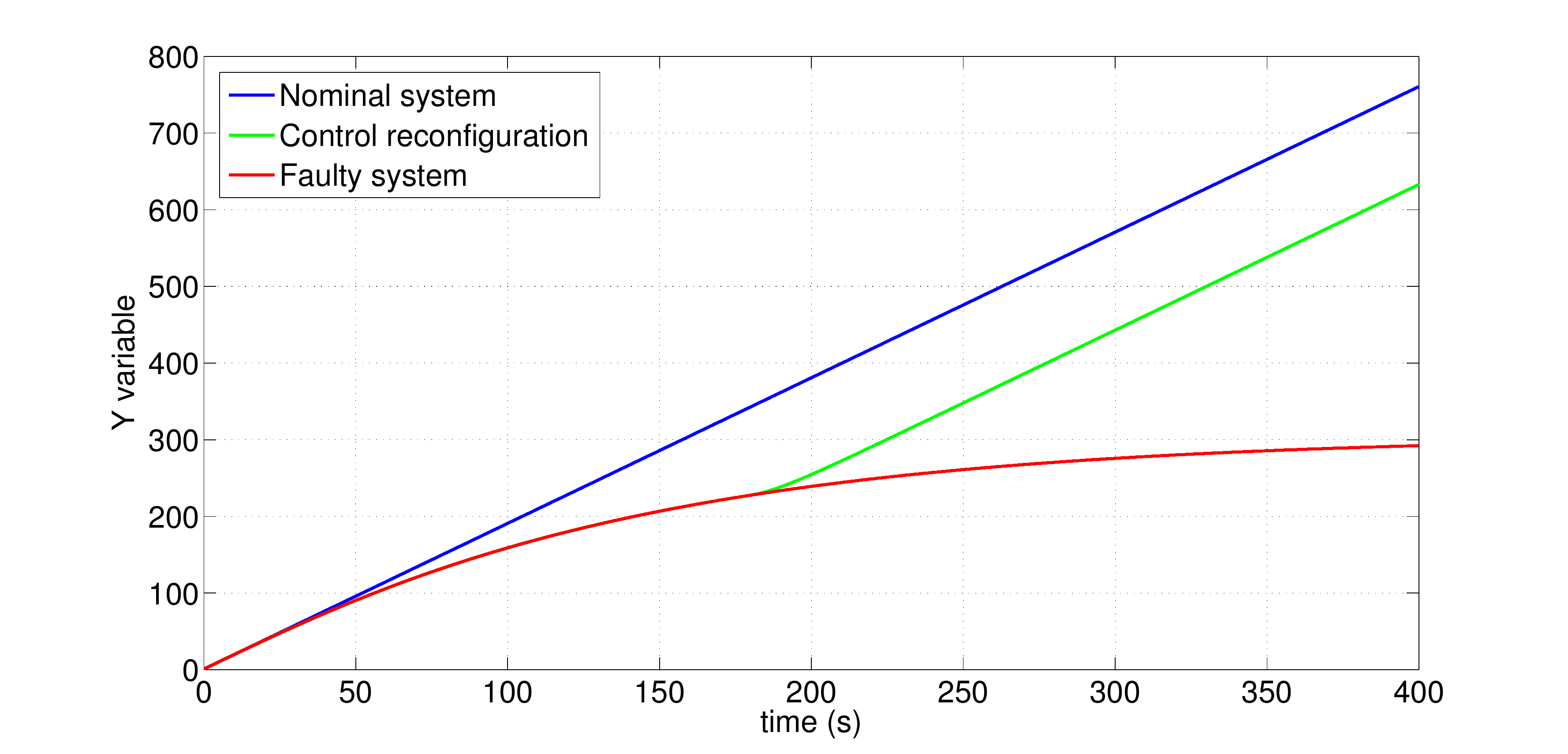}}
\caption{\label{fig6}Fault affecting the thruster $T_1$: evolution of the state variable $y_G(t)$}
\end{figure}

The control re-allocation can be performed through the reduced order pseudo-inverse method (\ref{simple-realloc}). Figures \ref{fig5}-\ref{fig6} show the evolution of the ship position in the case of a fault affecting the thruster $T_1$; the reconfiguration procedure is supposed to be activated after $t_0=180s$. The control is successfully re-allocated, and the commanded control effect is modified in order to track the original vessel velocities in surge and sway, while the rotational speed is automatically updated by the heading PID controller given by $a_\psi$. 

\end{document}